\documentclass[12pt, a4paper]{article}
\usepackage[margin=2.54cm]{geometry}

\usepackage[russian,english]{babel}
\usepackage{amsmath}
\usepackage[standard,thmmarks,amsmath]{ntheorem}
\usepackage{amssymb}
\usepackage{graphicx}
\usepackage[normalem]{ulem}
\usepackage{multirow, cite}

\newtheorem{note}{Note}
\newtheorem{conjecture}{Conjecture}
\def\le{\leqslant}

\def\ge{\geqslant}

\DeclareMathOperator*{\Argmin}{Arg\,min}

\title{\vspace*{2.5cm}\LARGE\bf Simple Alcohols with the Lowest Normal\\[2mm] Boiling Point Using
Topological Indices\thanks{This research is supported by the grant of Russian
Foundation for Basic Research, project No 13-07-00389.} }\author{\bf Mikhail Goubko,
\bf Oleg Miloserdov\\
\emph{\normalsize Institute of Control Sciences}\\
\emph{\normalsize of Russian Academy of Sciences, Moscow, Russia}\\
\texttt{\normalsize mgoubko@mail.ru, omilos92@gmail.com}}

\date{\vspace{0.1in}\normalsize (Received October 28, 2014)}

\begin{document}

\maketitle

\thispagestyle{empty}

\begin{abstract}
We find simple saturated alcohols with the given number of carbon atoms and the minimal normal boiling point.
The boiling point is predicted with a weighted sum of the generalized first Zagreb index, the second
Zagreb index, the Wiener index for vertex-weighted graphs, and a simple index caring for the degree of a carbon atom being
incident to the hydroxyl group. To find extremal alcohol molecules we characterize chemical trees of order $n$, which minimize
the sum of the second Zagreb index and the generalized
first Zagreb index, and also build chemical trees, which minimize the Wiener index over all chemical trees with given vertex weights.
\end{abstract}

\baselineskip=0.30in

\section{Introduction}

Consider a collection $\Omega$  of \emph{admissible molecules} (for example, represented with their structural formulas
or chemical graphs), each endowed with $k+1$ significant physical or chemical \emph{properties} (e.g., normal density,
normal boiling point, refraction coefficient, retention index, or more exotic and problem-specific ones), and let
$P_i(G)$, $i=0,...,k$, be the numeric value of the $i$-th property of a molecule $G\in \Omega$ (e.g., the normal
boiling point value). A typical problem of \emph{molecular design} is the following \emph{optimization problem}:
\begin{gather}\label{eq_property_optimization}
P_0(G)\rightarrow \min_{G\in \Omega}\,(\max_{G\in \Omega})\\
P_i^{min}\le P_i(G)\le P_i^{max}, i=1,...,k.\nonumber
\end{gather}

When the functions $P_i(\cdot)$ are only partially known from the
experiment, they are replaced with \emph{predicted} figures,
relating a chemical graph $G\in\Omega$ to the predicted value
$\tilde{P}_i(G)$ of the $i$-th physical or chemical property
($i=0,...,k$) by virtue of numeric characteristics (known as
\emph{molecular descriptors}), which can be calculated on basis of a
molecular structure. A typical \emph{quantitative structure-property
relation} (QSPR) includes several molecular descriptors, and is
presented as
$$\tilde{P}_i(G) = \tilde{P}_i(I_1(G),...,I_m(G)),i=0,...,k,$$
where $I_1(G),..., I_m(G)$ are the values of molecular descriptors for the molecular graph $G$. The simplest
\emph{linear regression} is just a weighted sum of descriptors:
$$\tilde{P}_i(G) = \alpha_{1,i} I_1(G)+...+\alpha_{m,i} I_m(G), i=0,...,k.$$

During recent decades a number of \emph{topological},
\emph{geometrical}, and \emph{quantum-mechanical} molecular
descriptors were suggested and studied \cite{Gutman10.1,Gutman10.2,
Todeschini00, Todeschini09}. Below we limit ourselves to topological
descriptors only (see, for instance, the handbook \cite{Balaban00})
to study problem (\ref{eq_property_optimization}) as a problem of
the \emph{extremal graph theory} \cite{Bollobas04} .

Exhaustive enumeration of all feasible molecules (the \emph{brute
force }approach) can only be used to solve this problem when the
feasible set is relatively small; for bigger sets mathematical
chemistry suggests a variety of limited search techniques. In
numerous papers lower and upper bounds of dozens topological indices
over various feasible sets were obtained \cite{Dobrynin01,
Volkmann02, Furtula13, Lin13, SurveyOptimLi08, ReviewOptimRada14,
Wang08, ReviewOptXu14, Zhang08, SurveyOptimZhou04}, and, in many
cases, extremal graphs were characterized. At the same time, the
problem (\ref{eq_property_optimization}) of optimization of a
composition of indices is still understudied.

In fact, finding lower and upper bounds of individual indices can be
a step towards solving problem (\ref{eq_property_optimization}), as
a linear combination of lower bounds is a lower-bound estimate of
the combination of indices. This estimate can be used in a
branch-and-bound algorithm of limited search. Yet, the quality of
the estimate may be considerably poor, resulting in lack of
efficient cuts in a branch-and-bound algorithm.

Anyway, the common shortcoming of an algorithmic approach to index
optimization is that it does not support the analysis of general
characteristics of an extremal molecule (i.e., of a corresponding
graph). When available, side information would be of great value on
why a certain graph is optimal or not, what shape the extremal
graphs have, etc. Such information is revealed using analytical
tools of discrete optimization.

In this paper we apply recent results in optimization of degree- and
distance-based topological indices to find a simple saturated
alcohol with the given molecular weight and minimal boiling point.
We reduce the property minimization problem to that of minimization
of a weighted linear combination of the generalized first Zagreb
index, the second Zagreb index, the vertex-weighted Wiener index,
and a simple index caring for the degree of a carbon atom being
incident to the hydroxyl group. Then we characterize minimizers of
this linear combination of indices (see
Fig.~\ref{fig_BP0_minimizers}) and of a couple of simpler
regressions (see Fig. \ref{fig_C_minimizers} and
\ref{fig_VWWI_minimizers}).

\section{Predicting Boiling Points of Simple Alcohols}\label{section_Huffman_alg}

The normal boiling point of a liquid is determined by its solvation
free energy. The solvation free energy can be predicted with high
accuracy from computer simulations (see \cite{Bren06, Bren12} for
details). The simulation-based approach solves well the ``direct
problem'' of predicting the solvation free energy for a given
molecule, but it does not help solving the ``inverse problem'' of
finding the molecule having the minimal solvation free energy (and,
consequently, the boiling point). For this reason we predict boiling
points of simple saturated alcohols (those having a general formula
$\textrm{C}_n\textrm{H}_{2n+1}\textrm{OH}$) with the aid of
topological indices.

Alcohols have relatively high boiling points when compared to
related compounds due to hydrogen bonds involving a highly polarized
hydroxyl group, and branched isomers have lower boiling points than
alcohols with the linear structure. Another structural feature
affecting the boiling point is the \emph{``oxygen shielding''
effect} \cite{AlcoholPenchev07}, when atoms surrounding the hydroxyl
group partially shield it preventing formation of hydrogen bonds
between molecules and, thus, decreasing the boiling point.

We considered several degree-based topological indices (the first Zagreb index $M_1$ \cite{GutmanTrin72}, the second
Zagreb index $M_2$ \cite{GutmanTrin72}, Randi\'c index \cite{Randic75} and the others), which are known to be good
metrics of branchiness, and, finally, the generalized first Zagreb index $C_1$ (see also \cite{Goubko14}) has shown the best results:
\begin{equation}
C_1(G)=\sum_{v\in V(G)}c(d_G(v)),
\end{equation}
where $V(G)$ is the vertex set of graph $G$, $d_G(v)$ is the degree of the vertex $v\in V(G)$ in graph $G$, and $c(d)$
is a non-negative function defined for degrees from 1 to 4. It can alternatively be written as
\begin{equation}\label{eq_gen1stZI}
C_1(G)=c(1) n_1(G)+c(2) n_2(G)+c(3) n_3(G)+c(4) n_4(G),
\end{equation}
where $n_i(G)$ is the atoms' count of degree $i=1,...,4$ in a molecular graph $G$, and $c(1), ..., c(4)$ are regression
parameters.

We also employed the classical \emph{second Zagreb index} \cite{GutmanTrin72}
\begin{equation}
M_2(G)=\sum_{uv\in E(G)}d_G(u)d_G(v),
\end{equation}
where $E(G)$ is the edge set of graph $G$.

Another index used was the \emph{Wiener index}, which had been the first topological index for boiling point prediction
\cite{Wiener47} due to its high correlation with the molecule's surface area. To account for heterogeneity of atoms we
allow each pair of vertices $u,v \in V(G)$ to have unique weight $\mu_G(u,v)$ and calculate the \emph{pair-weighted
Wiener index} as
\begin{equation}\label{eq_PWWI}
PWWI(G):=\frac{1}{2}\sum_{u,v\in V(G)}\mu_G(u,v) d_G(u,v),
\end{equation}
where $d_G(u, v)$ is the distance (the length of the shortest path) between vertices $u$ and $v$ in $G$. For example,
we can assign different weights to distances between pairs of carbon atoms and between carbons and the oxygen atom in an alcohol
molecule.

Regression tuning has shown the distances between carbon atoms to be
irrelevant for the alcohol boiling point, and only distances to the
oxygen matter. All such distances are accounted with equal weight,
so the pair-weighted Wiener index reduces to the \emph{distance} of
the oxygen atom, which was first used for the alcohol boiling point
prediction in \cite{AlcoholPenchev07}:
\begin{equation}\label{eq_WIO}
WI_\textrm{O}(G):=\sum_{u\in V(G)}d_G(u,\textrm{O}).
\end{equation}

In \cite{AlcoholPenchev07} a geometrical descriptor has also been suggested to account for oxygen shielding, but we
extend the approach by \cite{AlcoholJanezic06} instead, and introduce a simple topological index $S_i(G)$, which is
equal to unity when the carbon atom incident to the hydroxyl group in the alcohol molecule $G$ has degree $i=2,3,4$ (we
exclude \emph{methanol} from consideration), and is equal to zero otherwise.

We collected a data set of experimental boiling points under normal conditions for 79 simple saturated alcohols having
from 2 to 12 carbon atoms and representing various branchiness. Several data sources \cite{AlcoholGarcia08,
AlcoholJanezic06, KompanyZareh03, AlcoholPenchev07} were combined with priority on Alpha Aesar experimental data to
resolve discrepancy. In Table \ref{tab_stat_data} we present basic statistics about the data sample.

\begin{table}
\caption{Data sample: basic statistics} \label{tab_stat_data}
\begin{center}
\begin{tabular}{||c|c|c|c|c|c||}
  \hline
  \# carbons & \# isomers & min chain len. & max chain len. & min b.\,p. & max b.\,p.\\
  \hline
  2 & 1 & 2 & 2 & 78 & 78 \\
    \hline
  3 & 2 & 2 & 3 & 82.5 & 97 \\
    \hline
  4 & 4 & 2 & 4 & 82.4 & 117.5 \\
    \hline
  5 & 8 & 3 & 5 & 102 & 137 \\
    \hline
  6 & 17 & 3 & 6 & 120 & 157 \\
    \hline
  7 & 18 & 4 & 7 & 131.5 & 175.5\\
    \hline
  8 & 10 & 4 & 8 & 147.5 &  194\\
    \hline
  9 & 11 & 4 & 9 & 169.5 & 215\\
    \hline
  10 & 5 & 5 & 10 & 168 & 231\\
    \hline
  11 & 2 & 10 & 11 & 228.5 &  243\\
    \hline
  12 & 1 & 12 & 12 & 259 & 259\\
    \hline
  TOTAL & 79 & 2 & 12  & 78 & 259\\
  \hline
\end{tabular}
\end{center}
\end{table}

Information on boiling points of alcohols including more than 12 carbon atoms is less common and reliable. The complete
data set together with the best regressions is available online at \cite{GoubkoMilosOnline14}.

We randomly split the sample into the training set containing 50
cases and the testing set containing 29 cases. Then we examined
different linear regressions involving the descriptors mentioned
above.\footnote{ChemAxon Instant JChem$^\copyright$ was used for
index calculation. Authors would like to thank ChemAxon$^\circledR$
Ltd (http://www.chemaxon.com) for the academic license.} The best
performance and predictive power was obtained for the linear
combination of the oxygen's distance cube root
$WI_\textrm{O}(G)^\frac{1}{3}$ (with weight $b_1^0$), the
generalized first Zagreb index $C_1(G)$ (with weights
$c^0(1),...,c^0(4)$), the second Zagreb index $M_2(G)$ (weighted by
$b^0_3$), and the simple indicator of the sub-root's degree $S_2(G)$
(weighted by $b_2^0$):
\begin{equation}\label{eq_bp0}
BP^0=b_0^0+b_1^0 WI_\textrm{O}(G)^\frac{1}{3}+c^0(2) n_2(G)+c^0(3)
n_3(G)+b^0_2 S_2(G)+b^0_3 M_2(G).
\end{equation}

Below this regression is referred to as the \emph{basic} one. The
optimal values of weights $b_1^0, c^0(1),...,c^0(4), b_2^0, b_3^0$
(including the constant term $b_0^0$) calculated with the least
squares method are presented in the first column of Table
\ref{tab_param}. Parameters $c(1)$ and $c(4)$, which weight
variables $n_1(G)$ and $n_4(G)$ respectively, appear to be
insignificant and can be set to zero in (\ref{eq_gen1stZI}) when
calculating the generalized first Zagreb index. The precision of the
basic regression is shown in Table \ref{tab_precision}. It is
comparable to the best known relations
\cite{AlcoholGarcia08,AlcoholJanezic06,KompanyZareh03,
AlcoholPenchev07}. In the following sections we show how far we can
come in analytical and numeric minimization of this combination of
indices.

\begin{table}
\caption{Precision of regressions: Correlation coefficient and standard deviation (SD)}\label{tab_precision}
\begin{center}
\begin{tabular}{||c||c|c||c|c||c|c||}
  \hline
  \multirow{2}{*}{Data set} & \multicolumn{2}{c||}{Basic regression (``$^0$'')} & \multicolumn{2}{c||}{Regression I (``$^\textrm{I}$'')} & \multicolumn{2}{c||}{Regression II (``$^{\textrm{II}}$'')}\\
  \cline{2-7}
   & Corr. & SD, $^\circ\textrm{C}$ & Corr. & SD, $^\circ\textrm{C}$ & Corr. & SD, $^\circ\textrm{C}$\\
  \hline
  Training set & $0.997$ & $2.98$ & $0.997$ & $3.12$ & $0.996$ & $3.26$\\
  \hline
  Testing set & $0.996$ & $3.23$ & $0.995$ & $3.58$ & $0.994$ & $3.99$\\
  \hline
\end{tabular}
\end{center}
\end{table}

\begin{table}
\caption{Parameters of regressions} \label{tab_param}
\begin{center}
\begin{tabular}{|c||c|c|c||}
  \hline
  Coefficient & Basic regression (``$^0$'') & Regression I (``$^\textrm{I}$'') & Regression II (``$^{\textrm{II}}$'')\\
  \hline
  $b_0$ & $35.245$ & $50.626$ & $44.134$\\
  \hline
  $b_1$ & $12.233$ & $-$ & $3.851$\\
  \hline
  $b_2$ & $9.170$ & $11.295$ & $10.980$\\
  \hline
  $b_3$ & $1.486$ & $1.000$ & $-$\\
  \hline
  $c(1)$ & $-$ & $-$ & $-$\\
  \hline
  $c(2)$ & $9.514$ & $14.534$ & $17.727$ \\
  \hline
  $c(3)$ & $9.380$ & $20.172$ & $29.673$\\
  \hline
  $c(4)$ & $-$ & $17.015$ & $36.470$\\
  \hline
\end{tabular}
\end{center}
\end{table}

We also considered two simplifications of regression (\ref{eq_bp0}), for which alcohol molecules having minimal
predicted boiling point can be characterized analytically. The first one (below referred to as ``Regression I'') is
obtained by withdrawing $WI_\textrm{O}(G)$ in (\ref{eq_bp0}):
\begin{equation}\label{eq_bpI}
BP^\textrm{I}=b_0^\textrm{I}+c^\textrm{I}(2) n_2(G)+c^\textrm{I}(3) n_3(G)+c^\textrm{I}(4) n_4(G)+b^\textrm{I}_2
S_2(G)+b^\textrm{I}_3 M_2(G).
\end{equation}

See Table \ref{tab_precision} for the precision figures of
regression (\ref{eq_bpI}) under the values of parameters delivering
the best approximation to the training set (see the middle column of
Table \ref{tab_param}).

The second simplified regression (referred to as ``Regression II'' below) is obtained by withdrawing $M_2(G)$ in
(\ref{eq_bp0}):
\begin{equation}\label{eq_bpII}
BP^\textrm{II}=b_0^\textrm{II}+b_1^\textrm{II}
WI_\textrm{O}(G)^{\frac{1}{3}}+c^\textrm{II}(2)
n_2(G)+c^\textrm{II}(3) n_3(G)+c^\textrm{II}(4)
n_4(G)+b^\textrm{II}_2 S_2(G).
\end{equation}

In Table \ref{tab_precision} we show its precision under the optimal values of parameters depicted in the last
column of Table \ref{tab_param}.

The shortcoming of Regression II is that the term
$WI_\textrm{O}(G)^{\frac{1}{3}}$ appears to be insignificant after
disposal of $M_2(G)$, being responsible of approximately 1 per cent
of the residual sum of squares. Nevertheless, we keep this
regression for illustration of joint optimization of $C_1(G)$ and of
the pair-weighted Wiener index.

\section{Minimization of indices and their combinations}\label{section_min_comb}

In the present paper we find a simple saturated alcohol isomer with
$n-1$ carbon atoms having the lowest predicted boiling point. As the
regressions introduced in the previous section are tested only for
alcohols containing from 2 to 12 carbon atoms, we restrict our
attention to $n\le14$, where we can expect some accuracy of the
obtained results.

For $n\le 14$ admissible sets of all simple saturated alcohol molecules with $n-1$ carbons are not too extensive, and
allow for the brute-force enumeration. Moreover, we are sure that no aid of a computer is needed for an organic
chemist to draw a molecule being a good approximation to the boiling point minimizer for all $n\le 14$. However, our
aim is to show how analytic optimization techniques formalize the professional intuition and help making general
conclusions of verifiable reliability.

Let us characterize chemical trees minimizing indices introduced in the previous section and their combinations.

\subsection{Degree-based indices}\label{subsection_degree}

For a simple connected undirected graph $G$ denote with $W(G)$ the set of \emph{pendent vertices} (those having
degree 1) of the graph $G$, and with $M(G):=V(G)\backslash W(G)$ the set of \emph{internal vertices} (with degree $>1$) of $G$.

\begin{definition} A simple connected undirected graph of order $n$ is called a \emph{chemical tree} if it has $n-1$ edges and its vertex
degrees do not exceed $4$. Denote with $\mathcal{T}(n)$ the set of all chemical trees of order $n$.
\end{definition}

\begin{definition} A \emph{pendent-rooted chemical tree} is a chemical tree, in which one pendent vertex is
distinguished and called a \emph{root}. A typical pendent-rooted tree is denoted with $T_r$, with $r$ being its root. A
vertex being incident to the root in $T_r$ is called a \emph{sub-root} and is denoted as $sub(T_r)$. Denote with
$\mathcal{R}(n)$ the set of all pendent-rooted chemical trees of order $n$.
\end{definition}

Clearly, the set $\Omega(n-1)$ of all molecules of simple saturated alcohols having $n-1$ carbons coincides with the set
$\mathcal{R}(n)$ of pendent-rooted chemical trees of order $n$ (with a root corresponding to the hydroxyl group and
the other vertices forming the carbon skeleton of a molecule).

For a topological index $I(\cdot)$ defined on an admissible set of graphs $\mathcal{G}$ introduce the notation
$I_\mathcal{G}^*:=\min_{G\in \mathcal{G}}I(G)$ and let $\mathcal{G}_I^*:=\Argmin_{G\in \mathcal{G}}I(G)$
be the set of graphs minimizing $I(\cdot)$ over $\mathcal{G}$. For example, $\mathcal{T}_{M_2}^*(n)$ is the set of chemical
trees of order $n$ minimizing the second Zagreb index $M_2(\cdot)$.

Define also the set $\mathcal{R}_i(n):=\{T\in \mathcal{R}(n): d_T(sub(T))=i\}$ of all pendent-rooted trees with
a sub-root having degree $i=2,...,4$.

We start with the following obvious statement.

\begin{lemma}\label{lemma_S_i}$S_i(G)$ achieves its minimum at any pendent-rooted chemical tree with sub-root's degree other than $i$. In other
words, $\mathcal{R}_{S_i}^*(n)=\mathcal{R}(n)\backslash \mathcal{R}_i(n)$.
\begin{proof} is straightforward, as $S_i(G)=1$ for all $G\in \mathcal{R}_i(n)$, and $S_i(G)=0$ otherwise.
\end{proof}
\end{lemma}

Indices $C_1(G)$ and $M_2(G)$ do not account for heterogeneity of
atoms in a molecule, so we can minimize them over the set
$\mathcal{T}(n)$ of all chemical trees of order $n$ and then assign
the root to an arbitrary pendent vertex of the index-minimizing tree
to obtain a pendent-rooted tree, which minimizes the index.

Consider an ``ad-hoc'' degree-based topological index
\begin{equation}\label{eq_ad_hoc_index}
C(G):=C_1(G)+b_3M_2(G)=\sum_{v\in V(G)}c(d_G(v))+b_3\sum_{uv\in
E(G)}d_G(u)d_G(v),
\end{equation}
where $b_3$ is an arbitrary real constant (we keep notation $b_3$
for compatibility with equations (\ref{eq_bp0}), (\ref{eq_bpI})).

\begin{definition}
A chemical tree $T\in \mathcal{T}(n)$ is \emph{extremely branched},
if its internal vertices have degree $4$, except one vertex having
degree $2$ when $n\bmod{3}=0$, or one vertex having degree $3$ when
$n\bmod{3}=1$.
\end{definition}

\begin{theorem}\label{theorem_C_sufficient}
Assume the following inequalities hold:
\begin{gather}
c(1)+c(4)+18b_3<c(2)+c(3),\label{eq_degree_index_condition23}\\
c(1)+c(3)+8b_3<2c(2),\label{eq_degree_index_condition22}\\
c(2)+c(4)+8b_3<2c(3).\label{eq_degree_index_condition33}
\end{gather}
If a chemical tree $T\in \mathcal{T}(n)$ for $n\ge 3$ minimizes $C(\cdot)$ over all chemical trees from $\mathcal{T}(n)$,
then $T$ is an extremely branched tree.
For $n\le 17$ the inequality \emph{(\ref{eq_degree_index_condition23})} can be weakened to
\begin{equation}
c(1)+c(4)+17b_3<c(2)+c(3).\label{eq_degree_index_condition23bis}
\end{equation}
\begin{proof}
We employ the standard argument of index monotonicity with respect to certain tree transformations.
Assume the theorem does not hold, and vertices $u,v \in M(T)$ exist such that $u\neq v$ and $d_T(u),d_T(v)<4$. Four
cases are possible.
\begin{enumerate}
    \item $d_T(u)=d_T(v)=2$. Let $v_1,v_2\in V(T)$ be the vertices incident to $v$ in $T$.
    Also, let $u_1,u_2\in V(T)$ be the vertices incident to $u$ in $T$, and $u_2$ lies on the path to the vertex $v$ in $T$. Without loss of generality assume that
    \begin{equation}\label{eq_case1_ineq}
    d_T(u_1)+d_T(u_2)\ge d_T(v_1)+d_T(v_2).
    \end{equation}
    Consider a
    graph $T'\in \mathcal{T}(n)$ obtained from $T$ by replacing the edge $u_1u$ with the edge $u_1v$. It is easy to see
    that $T'$ is a tree.
    The degree of the vertex $u$ in $T'$ is decreased by one, the degree of vertex $v$ is increased by one, therefore, if $u_2\neq v$, we have
    $$C(T')-C(T)=c(1)+c(3)-2c(2)+b_3(d_T(u_1)+d_T(v_1)+d_T(v_2)-d_T(u_2)).$$
    From (\ref{eq_case1_ineq}) we obtain $C(T')-C(T)\le c(1)+c(3)-2c(2)+2b_3d_T(u_1)$. Since
    vertex degrees $\le4$ in a chemical tree, from (\ref{eq_degree_index_condition22}) we have $C(T')-C(T)\le c(1)+c(3)-2c(2)+8b_3<0$,
    which contradicts the assumption that $T$ minimizes $C(\cdot)$.

    If $u_2=v$, in the same manner obtain $C(T')-C(T)\le c(1)+c(3)-2c(2)+7b_3$. From (\ref{eq_degree_index_condition22}), it is also negative.
    \item $d_T(u)=2$, $d_T(v)=3$. Let $v_1,v_2,v_3\in V(T)$ be the vertices incident to $v$ in $T$.
    Also, let $u_1,u_2\in V(T)$ be the vertices incident to $u$ in $T$, and assume $u_2$ lies on the path to the vertex $v$ in the tree $T$. Consider a
    tree $T'\in \mathcal{T}(n)$ obtained from $T$ by replacing the edge $u_1u$ with the edge $u_1v$.
    By analogy to the previous case, if $u_2\neq v$, we obtain $C(T')-C(T)=c(1)+c(4)-c(2)-c(3)+b_3(2d_T(u_1)+d_T(v_1)+d_T(v_2)+d_T(v_3)-d_T(u_2)).$
    Vertex degrees $\le4$ in a chemical tree. Moreover, $d_T(u_2)\ge 2$, since it is an intermediate vertex on
    the path $u,u_2,...,v$. Therefore,
    $$C(T')-C(T)\le c(1)+c(4)-c(2)-c(3)+18b_3,$$
    and, from (\ref{eq_degree_index_condition23}), $C(T')-C(T)<0$, which is a contradiction.

    To prove the weaker inequality (\ref{eq_degree_index_condition23bis}) we are enough to prove that
    $C(T')-C(T)<0$ for $n\le 17$, since $d_T(u_1)=d_T(v_1)=d_T(v_2)=d_T(v_3)=4$,
    and $d_T(u_2)=2$ is possible only in a tree of order 18 or more
    (an example is depicted in Fig.~\ref{fig_tree18}a), and $C(T')-C(T)= c(1)+c(4)-c(2)-c(3)+17b_3$
    for the tree $T$ depicted in Fig.~\ref{fig_tree18}b.
\begin{figure}[htpb]
\begin{center}
  \includegraphics[width=10cm]{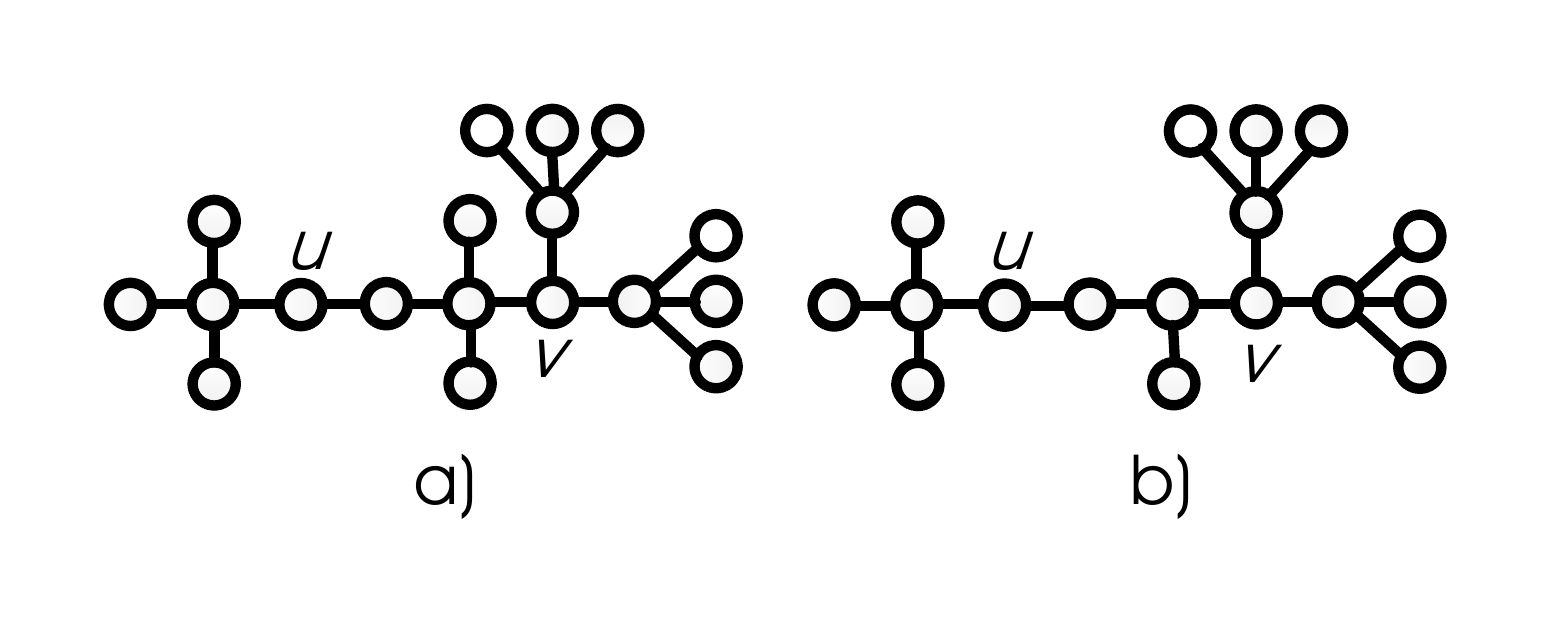}
  \caption{To the proof of inequality (\ref{eq_degree_index_condition23bis})}
  \label{fig_tree18}
\end{center}
\end{figure}
    If $u_2=v$, without loss of generality assume that $u=v_3$.
    Then $C(T')-C(T)=c(1)+c(4)-c(2)-c(3)+b_3(2d_T(u_1)+d_T(v_1)+d_T(v_2)-2)\le c(1)+c(4)-c(2)-c(3)+14b_3$.
    From (\ref{eq_degree_index_condition23}) (or from (\ref{eq_degree_index_condition23bis}), if $n\le 17$), it is negative.
    \item The case of $d_T(v)=2$, $d_T(u)=3$ is considered in the same manner.
    \item $d_T(u)=d_T(v)=3$. Let $v_1,v_2,v_3\in V(T)$ be the vertices incident to $v$ in $T$.
    Also, let $u_1,u_2,u_3\in V(T)$ be the vertices incident to $u$ in $T$, with $u_1$ not laying on a path to $v$ in $T$.
    Without loss of generality assume that
    \begin{equation}\label{eq_case4_ineq}
    d_T(u_1)+d_T(u_2)+d_T(u_3)\ge d_T(v_1)+d_T(v_2)+d_T(v_3).
    \end{equation}
    Consider a tree $T'\in \mathcal{T}(n)$ obtained from $T$ by replacing the edge $u_1u$ with the edge $u_1v$.
    If $uv\notin E(T)$, then from (\ref{eq_case4_ineq}) and $d_T(u_1)\le 4$ we have
    $$C(T')-C(T)=c(2)+c(4)-2c(3)+b_3(d_T(u_1)+d_T(v_1)+d_T(v_2)+d_T(v_3)-d_T(u_2)-$$
    $$-d_T(u_3))\le c(2)+c(4)-2c(3)+2b_3d_T(u_1)\le c(2)+c(4)-2c(3)+8b_3,$$
    which is less than zero due to (\ref{eq_degree_index_condition33}), and $T$ cannot minimize $C(\cdot)$.
    If $uv\in E(T)$, in the same way deduce $C(T')-C(T)\le c(2)+c(4)-2c(3)+7b_3$, which is negative.
\end{enumerate}

The obtained contradictions prove that no more than one internal vertex in $T$ may have degree less than 4.

As $T\in \mathcal{T}(n)$ and $n>1$, the well-known equity holds:
\begin{equation}\label{eq_tree_degrees}
n_1(T)+2n_2(T)+3n_3(T)+4n_4(T)=2(n-1).
\end{equation}
On the other hand,
\begin{equation}\label{eq_tree_vertices}
n_1(T)+n_2(T)+n_3(T)+n_4(T)=n.
\end{equation}

Assume that $n_2(T)=1$, so that $n_3(T)=0$. From (\ref{eq_tree_vertices}) we have $n_1(T)+n_4(T)=n-1$, therefore,
(\ref{eq_tree_degrees}) makes $n=3+3n_4(T)$ and, since $n_4(T)\in \mathbb{N}_0$, $n\bmod{3}=0$.

In the same manner we show that if $n_3(T)=1$ then $n\bmod{3}=1$. If both $n_2(T)$ and $n_3(T)=0$, then $n\bmod{3}=2$, and
the proof is complete.
\end{proof}
\end{theorem}

\begin{corollary}\label{corollary_C1}
Under conditions of Theorem \emph{\ref{theorem_C_sufficient}}, any tree $T$ minimizing
$C(T)=C_1(T)+b_3M_2(T)$ over $\mathcal{T}(n)$ enjoys the same number $n_i$ of vertices of degree $i=1,...,4$. Therefore, $C_1(T) = C_1(T')$
for any pair of trees $T, T'\in \mathcal{T}_C^*(n)$.
\end{corollary}

\begin{figure}[htpb]
\begin{center}
  \includegraphics[width=12cm]{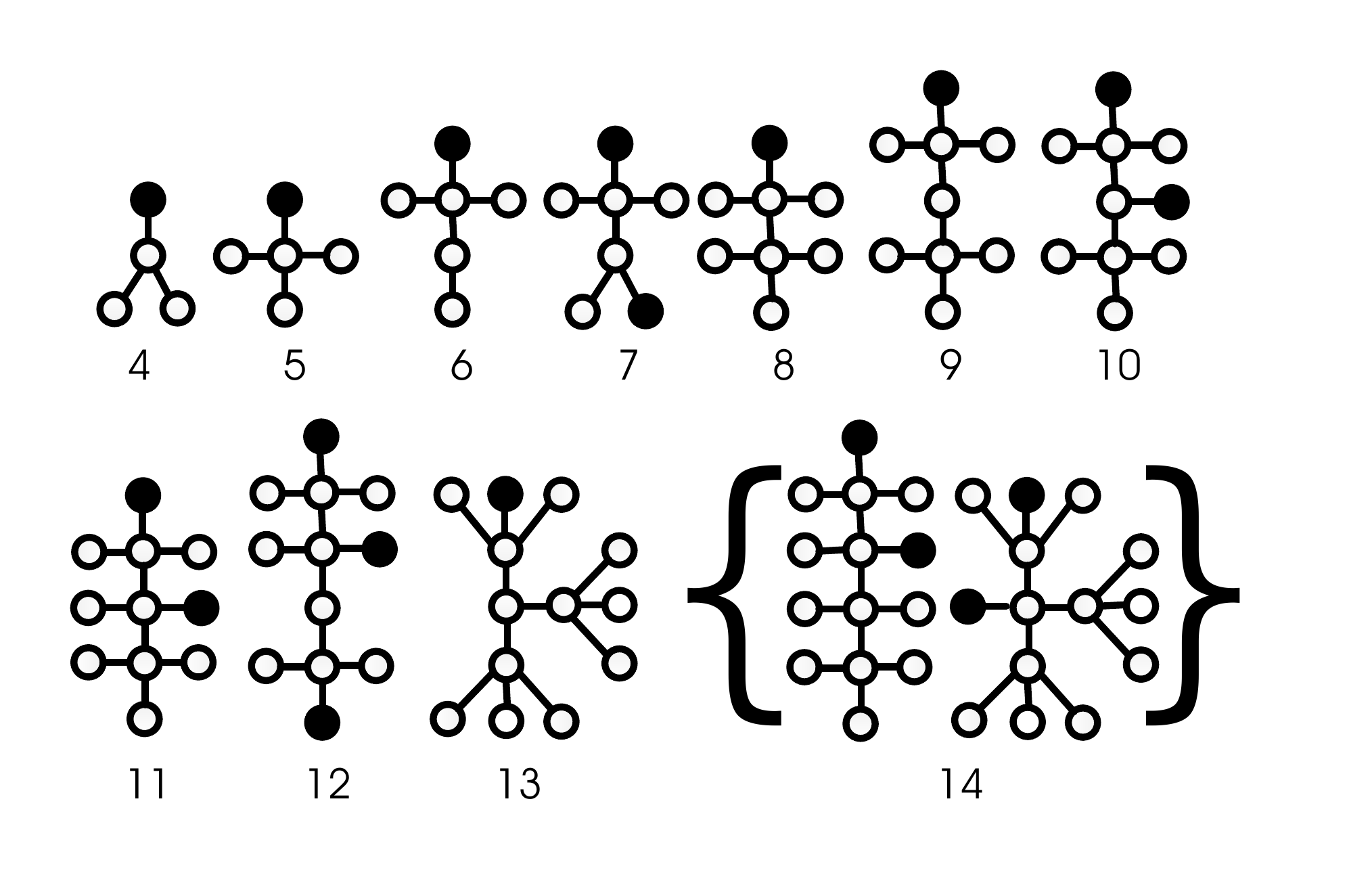}
  \caption{Chemical trees minimizing the ``ad-hoc'' index $C(\cdot)$ for $n=4,...,14$
  (alcohol molecules minimizing $BP^\textrm{I}(\cdot)$ with possible oxygen positions filled with black)}
  \label{fig_C_minimizers}
\end{center}
\end{figure}

\begin{corollary}\label{corollary_ad_hoc}
Under conditions of Theorem \emph{\ref{theorem_C_sufficient}} the sets $\mathcal{T}_C^*(n)$ for $n=4,...,14$ are
depicted in Fig. \emph{\ref{fig_C_minimizers}}. $\mathcal{T}_C^*(n)$ contains the sole tree for $n<14$ , while
$\mathcal{T}_C^*(14)$ contains two trees.
\begin{proof}
From Corollary~\ref{corollary_C1} we learn that only the value of $M_2(\cdot)$ may vary within
$\mathcal{T}_C^*(n)$.

From Theorem \ref{theorem_C_sufficient}, for $n\in \{5,8,11,14\}$ an optimal tree is a 4-tree (in which all internal
vertices have degree 4). Each of $n_1$ \emph{stem} edges (those incident to a pendent vertex) adds 4 to the value of
$M_2$, while each of $n_4-1$ edges connecting internal vertices adds 16 to the value of $M_2$. Since $n_1$ and $n_4$
are fixed for fixed $n$, all 4-trees have the same value of $M_2(\cdot)$ (and, therefore, the same value of
$C(\cdot)$). Consequently, for for $n = 5,8,11,14$ the set $\mathcal{T}_C^*(n)$ consists of all 4-trees of order $n$
(see Fig. \ref{fig_C_minimizers}).

If $T\in\mathcal{T}_C^*(n)$, and $n \in \{6,9,12\}$, one internal vertex $u\in M(T)$ has degree $d_T(u)=2$, while all
others have degree 4. For $n=6$ only one such tree exists depicted in Fig. \ref{fig_C_minimizers}. It is easy to
check that $M_2(\cdot)$ is minimized if vertex $u$ is incident to two internal vertices. Only
one such tree exists for $n=9$ (see Fig. \ref{fig_C_minimizers}), and the same is true for $n=12$.

For $n \in \{4,7,10,13\}$ any tree $T\in\mathcal{T}_C^*(n)$ has one internal vertex $u\in M(T)$ of degree $d_T(u)=3$,
while all other have degree 4. For $n=4$ only one such tree exists depicted in Fig. \ref{fig_C_minimizers}, and the
same is true for $n=7$. Again, it is easy to check that, in the context of $M_2(\cdot)$ minimization, vertex $u$ being
incident to three internal vertices is strictly preferable to vertex $u$ being incident to one pendent and two internal
vertices, which is, in turn, preferred to $u$ having two incident pendent vertices. So, optimal trees for $n=10,13$ are
depicted in Fig. \ref{fig_C_minimizers} (black and white filling of circles is explained below).
\end{proof}
\end{corollary}

The same logic allows continuing the sequence of $C(\cdot)$-minimizers to $n > 14$.

Therefore, Theorem \ref{theorem_C_sufficient} says that, when conditions (\ref{eq_degree_index_condition23})-(\ref{eq_degree_index_condition33}) hold,
chemical trees minimizing $C(\cdot)$ have as many vertices of maximal degree 4 as possible. A similar result can be
proved for the modified Wiener index $WI_\textrm{O}(\cdot)$.

\subsection{Wiener index}\label{section_wiener}

A simple connected undirected graph $G$ is called
\emph{vertex-weighted} if each vertex $v\in V(G)$ is endowed with a
non-negaive weight $\mu_G(v)$. With $\mu_G$ we denote the total
vertex weight of the graph $G$, and $\mathcal{WT}(n)$ stands for the
set of all vertex-weighted trees of order $n$.

Klav\v{z}ar and Gutman \cite{Klavzar97} defined the \emph{Wiener index for vertex-weighted graphs} as
$$VWWI(G):=\frac{1}{2}\sum_{u,v\in V(G)}\mu_G(u)\mu_G(v) d_G(u,v).$$

Clearly, $VWWI(\cdot)$ is a special case of the pair-weighted Wiener index $PWWI(\cdot)$ (defined with formula
(\ref{eq_PWWI})) for $\mu_G(u, v):=\mu_G(u)\mu_G(v)$. The path-weighted Wiener index is poorly studied at the moment,
but, fortunately, $WI_\textrm{O}(\cdot)$, which is the point of our current interest, can be reduced to the Wiener
index for vertex-weighted graphs.

For every alcohol molecule from $\Omega(n-1)$ (or, equivalently, for every pendent-rooted tree $T_r\in\mathcal{R}(n)$)
define a vertex-weighted tree $T(\varepsilon)\in \mathcal{WT}(n)$ by assigning the weight
$\mu_{T(\varepsilon)}(v):=\varepsilon$ to each vertex $v\in V(T_r)$ (a carbon atom) except the root $r$, and assigning
the weight $\mu_{T(\varepsilon)}(r):=1/\varepsilon$ to the root (the oxygen atom). It is easy to see that under these
weights $\lim_{\varepsilon\rightarrow 0} VWWI(T(\varepsilon))=WI_\textrm{O}(T_r)$. Since $WI_\textrm{O}(\cdot)$ is
integer-valued, minimizers of $VWWI(\cdot)$ and of $WI_\textrm{O}(\cdot)$ coincide for sufficiently small
$\varepsilon$.

In \cite{Goubko15} the majorization technique suggested by Zhang et
al. \cite{Zhang08} is used to minimize $VWWI(\cdot)$  over the set
of trees with given vertex weights and degrees. Below we recall the
notation and selected theorems from \cite{Goubko15}. We use them to
find the extremal vertex degrees over the set of all trees of order
$n$ with fixed vertex weights.

\begin{definition}
Consider a vertex set $V$. Let the function $\mu: V \rightarrow
\mathbb{R}_+$ assign a non-negative weight $\mu(v)$ to each vertex
$v\in V$, while the function $d: V \rightarrow \mathbb{N}$ assign a
natural degree $d(v)$. The tuple $\langle\mu, d\rangle$ is called
\emph{a generating tuple} if the following identity holds:
\begin{equation}\label{tree_degrees_identity}
\sum_{v\in V} d(v)=2(|V|-1).
\end{equation}

Denote with $\overline{\mu}:=\sum_{v\in V}\mu(v)$ the total weight of the vertex set $V$. Let
$\mathcal{WT}(\mu):=\{T\in \mathcal{WT}(|V|): V(T)=V, \mu_T(v)=\mu(v) \text{ for all }v\in V\}$ be the set of trees
over the vertex set $V$ with vertex weights $\mu(\cdot)$. For the set $\mathcal{WT}(\mu, d):=\{T\in \mathcal{WT}(\mu):
d_T(v)=d(v) \text{ for all }v\in V\}$ we also require vertices to have degrees $d(\cdot)$.
\end{definition}

Let $V(\mu, d)$ be the domain of functions of a generating tuple $\langle\mu, d\rangle$. Introduce the set $W(\mu,
d):=\{w \in V(\mu, d): d(w)=1\}$ of \emph{pendent} vertices and the set $M(\mu, d):=V(\mu,d)\backslash W(\mu,d)$ of
\emph{internal} vertices.

\begin{definition}
We will say that in a generating tuple $\langle\mu, d\rangle$
\emph{weights are degree-monotone}, if for any $m, m' \in M(\mu,d)$
from $d(m) < d(m')$ it follows that $\mu(m) \le \mu(m')$, and for
any $w\in W(\mu,d)$ we have $\mu(w)>0$.
\end{definition}


For a generating tuple $\langle\mu,d\rangle$ the \emph{generalized Huffman algorithm} \cite{Goubko15} builds a tree $H \in
\mathcal{WT}(\mu,d)$ as follows.

\textbf{Setup.} Define the vertex set $V_1 := V(\mu,d)$ and the functions $\mu^1$ and $d^1$, which endow its vertices
with weights $\mu^1(v) := \mu(v)$ and degrees $d^1(v) := d(v)$, $v \in V_1$. We start with the empty graph $H$ over the
vertex set $V(\mu,d)$.

\textbf{Steps $i = 1, ..., q-1$.} Denote with $m_i$ the vertex having the least degree among the vertices of the least weight in $M(\mu^i,d^i)$. Let $w_1,
..., w_{d(m_i)-1}$ be the vertices having $d(m_i)-1$ least weights in $W(\mu^i,d^i)$. Add to $H$ edges $w_1m_i, ...,
w_{d(m_i)-1}m_i$.

Define the set $V_{i+1} := V_i \backslash \{w_1, ..., w_{d(m_i)-1}\}$ and functions $\mu^{i+1}(\cdot), d^{i+1}(\cdot)$, endowing its
elements with weights and degrees as follows:
\begin{gather}
\mu^{i+1}(v) := \mu^i(v)\text{ for }v \neq m_i, \hspace{20pt} \mu^{i+1}(m_i) := \mu^i(m_i)+\mu^i(w_1)+...+\mu^i(w_{d(m_i)-1}),\nonumber\\
\label{eq_Huffman_tuples}d^{i+1}(v) := d^i(v)\text{ for }v \neq m_i, \hspace{20pt} d^{i+1}(m_i) := 1.
\end{gather}

\textbf{Step $q$.} Consider a vertex $m_q \in M(\mu^q, d^q)$. By construction, $|M(\mu^q, d^q)|=1$, $|W(\mu^q,
d^q)|=d(m_q)$. Add to $H$ edges connecting all vertices from $W(\mu^q, d^q)$ to $m_q$. Finally, set $\mu_H(v) :=
\mu(v)$, $v \in V(H)$.

\begin{theorem}\label{theorem_huffman_vwwi}\textbf{\emph{\cite{Goubko15}}} \sloppy If weights are degree-monotone
in a generating tuple $\langle\mu,d\rangle$, then $T\in \mathcal{WT}_{VWWI}^*(\mu,d)$ if and only if
$T\in\mathcal{WT}(\mu,d)$ and $T$ is a Huffman tree. In other words, only a Huffman tree minimizes the Wiener
index over the set of trees whose vertices have given weights and degrees.
\end{theorem}

In the present subsection we study how the value of
$VWWI_{\mathcal{WT}}^*(\mu,d)$ changes with degrees $d(\cdot)$. Our
results are analogous to those proved by Zhang et al. \cite{Zhang08}
for the ``classical'' Wiener index. Following \cite{Goubko15}, we
reformulate the problem for directed trees.

\begin{definition}
A (weighted) \emph{directed tree} is a weighted connected directed graph with each vertex except the \emph{terminal
vertex}\footnote{Typically it is called a \emph{root}, but we will use an alternative notation to avoid confusion with
a root of a pendent-rooted tree introduced in the previous subsection.} having the sole outbound arc and the terminal
vertex having no outbound arcs.
\end{definition}

An arbitrary tree $T \in \mathcal{WT}(n)$ can be transformed into a directed tree by choosing an internal vertex
$t\in M(T)$, and replacing all its edges with arcs directed towards (a terminal vertex) $t$. Let us denote with
$\mathcal{WD}$ the collection of all directed trees, which can be obtained in such a way, and let $\mathcal{WD}(\mu,
d)$ stand for all directed trees obtained from $\mathcal{WT}(\mu,d)$. Vice versa, in a directed tree from
$\mathcal{WD}(\mu,d)$ replacement of all arcs with edges makes some tree from $\mathcal{WT}(\mu, d)$.

If at Step $i=1,...,q$ of the generalized Huffman algorithm we add arcs towards the vertex $m_i$ (instead of undirected
edges), we obtain a \emph{directed Huffman tree} with the terminal vertex $m_q$.


\begin{definition}
For a vertex $v\in V(T)$ of a directed tree $T\in \mathcal{WD}$
define its \emph{subordinate group} $g_T(v)\subseteq V(T)$ as the
set of vertices having the directed path to the vertex $v$ in the
tree $T$ (the vertex $v$ itself belongs to $g_T(v)$). The
\emph{weight $f_T(v)$ of a subordinate group} $g_T(v)$ is defined as
the total vertex weight of the group: $f_T(v):=\sum_{u\in g_T(v)}
\mu_T(u)$.
\end{definition}

\begin{note}\label{note_positive}If all pendent vertices in $T$ have positive weights,
then $f_T(v)>0$  for any $v\in V(T)$. In particular, it is true for
any $T\in \mathcal{WD}(\mu, d)$, if weights in $\langle\mu,d\rangle$
are degree-monotone.
\end{note}

The Wiener index is defined for directed trees by analogy to the
case of undirected trees: we simply ignore the arcs' direction when
calculating distances. Therefore, a tree and a corresponding
directed tree share the same value of the Wiener index.

The value of the Wiener index for a directed tree $T_t\in
\mathcal{WD}(\mu)$ with a terminal vertex $t\in M(T)$ can be written
\cite{Goubko15} as:
\begin{equation}\label{eq_VWWI_directed}
VWWI(T_t)=\sum_{v\in V(T_t)\backslash
\{t\}}f_{T_t}(v)(\bar{\mu}-f_{T_t}(v))=\sum_{v\in V(T)\backslash
\{t\}}\chi(f_{T_t}(v)),
\end{equation}
where $\chi(x):=x(\bar{\mu}-x)$, and thus, the problems of Wiener index minimization for vertex-weighted trees and for
weighted directed trees are equivalent.

\begin{definition}\label{def_vector}
Every directed tree $T$ is associated with the \emph{vector of subordinate groups' weights} $\mathbf{f}(T):=(f_T(v))_{v\in
V(T)\backslash\{t\}}$, where $t$ is the terminal vertex of $T$. From
equation (\ref{eq_VWWI_directed}) we see that the vector $\mathbf{f}(T)$ completely  determines the value of $VWWI(T)$.
\end{definition}

\begin{definition}\textbf{ \cite{Marshall79,Zhang08}}
For the real vector $\mathbf{x}=(x_1,..., x_p)$, $p\in \mathbb{N}$,
denote with $\mathbf{x}_\uparrow=(x_{[1]},..., x_{[p]})$ the vector
where all components of $\mathbf{x}$ are arranged in ascending
order.
\end{definition}

\begin{definition} \textbf{\cite{Marshall79,Zhang08}}
A non-negative vector $\mathbf{x} = (x_1, ..., x_p)$, $p \in
\mathbb{N}$, \emph{weakly majorizes} a non-negative vector
$\mathbf{y} = (y_1, ..., y_p)$ (which is denoted with
$\mathbf{x}\succeq \mathbf{y}$) if
$$\sum_{i=1}^k x_{[i]} \le \sum_{i=1}^k y_{[i]} \text{ for all }k=1,...,p.$$
If $\mathbf{x}_\uparrow \neq \mathbf{y}_\uparrow$, then $\mathbf{x}$
is said to \emph{strictly weakly majorize} $\mathbf{y}$ (which is
denoted with $\mathbf{x}\succ \mathbf{y}$).
\end{definition}

We will need the following properties of weak majorization.

\begin{lemma}\label{lemma_Zhang_b}\textbf{\emph{\cite{Marshall79,Zhang08}}}
Consider a positive number $b>0$ and two non-negative vectors
$\mathbf{x} = (x_1, ..., x_k, y_1, ..., y_l)$ and $\mathbf{y} = (x_1
+ b, ..., x_k + b, y_1 - b, ..., y_l - b)$, such that $0 \le k \le
l$. If $x_i \ge y_i$ for $i = 1, ..., k$, then $\mathbf{x}\prec
\mathbf{y}$.
\end{lemma}

\begin{lemma}\label{lemma_Zhang_xy}\textbf{\emph{\cite{Marshall79,Zhang08}}}
If $\mathbf{x}\preceq \mathbf{y}$ and $\mathbf{x}' \prec
\mathbf{y'}$, then $(\mathbf{x},\mathbf{x'})\prec
(\mathbf{y},\mathbf{y'})$, where $(\mathbf{x},\mathbf{x'})$ means
concatenation of vectors $\mathbf{x}$ and $\mathbf{x'}$.
\end{lemma}

\begin{lemma}\label{lemma_Zhang_concave}\textbf{\emph{\cite{Marshall79,Zhang08}}}
If $\chi(x)$ is a increasing concave function, and
$(x_1,...,x_p)\succeq (y_1,...,y_p)$, then $\sum_{i=1}^p\chi(x_i)\le
\sum_{i=1}^p\chi(y_i)$, and equality is possible only when
$(x_1,...,x_p)_\uparrow\\=(y_1,...,y_p)_\uparrow$.
\end{lemma}

The following lemma establishes an important property of directed Huffman trees:

\begin{lemma}\label{lemma_Huffman_monotonicity}\textbf{\emph{\cite{Goubko15}}}
For any directed Huffman tree $H$
\begin{equation}\label{eq_weights_monotone}
vm,v'm'\in E(H), m \neq m', f_H(v)<f_H(v') \Rightarrow f_H(m) < f_H(m').
\end{equation}
\end{lemma}

\begin{lemma}\label{lemma_major}
Consider generating tuples $\langle\mu, d\rangle$ and $\langle\mu,
d'\rangle$ defined on the same vertex set, and let weights be
degree-monotone in $\langle\mu, d\rangle$. Let the values of degree
functions $d(\cdot)$ and $d'(\cdot)$ differ only for vertices $u$
and $v$, such that $d(u)\ge d(v)$ and $\mu(u)\ge \mu(v)$, while
$d'(u)=d(u)+1$, $d'(v)=d(v)-1$. Then for every directed tree $T\in
\mathcal{WD}(\mu,d)$ there exists such a directed tree $T'\in
\mathcal{WD}(\mu,d')$ that $\mathbf{f}(T')\succ \mathbf{f}(T)$.
\begin{proof}
By Theorem 2 from \cite{Goubko15}, such a directed Huffman tree
$H\in \mathcal{WD}(\mu,d)$ exists, that $\mathbf{f}(H)\succeq
\mathbf{f}(T)$.\footnote{Please note the different notation in
\cite{Goubko15} ($x\preceq_w y$ is used in \cite{Goubko15} where we
write $x\succeq y$).} Since $H\in \mathcal{WD}(\mu,d)$ and
$d'(v)=d(v)-1\ge1$, we know that $d(v)\ge 2$ and $v$ has an incoming
arc in $H$ from some vertex $v'\in V$. Weights are degree-monotone
in $\langle\mu,d\rangle$, therefore, by construction of a Huffman
tree, $f_H(u)\ge f_H(v)$ and, thus, without loss of generality we
can assume that $u \notin g_H(v)$.

Assume that $v \in g_H(u)$. Then a path $(v, m_1, ..., m_l, u)$
exists in $H$ from the vertex $v$ to the vertex $u$, where $l\ge0$.
Consider a directed tree $T'$ obtained from $H$ by deleting the arc
$v'v$ and adding the arc $v'u$ instead. It is clear that $T'\in
\mathcal{WD}(\mu,d')$, and weights of groups subordinated to
vertices $v, m_1, ..., m_l$ decrease by $f_H(v')$ (which is positive
by Note \ref{note_positive}), while weights of the other vertices do
not change. Therefore, by Lemma \ref{lemma_Zhang_b},
$$\mathbf{y} := (f_{T'}(v), f_{T'}(m_1), ..., f_{T'}(m_l))=$$
$$=(f_H(v) -f_H(v'), f_H (m_1) -f_H(v'), ..., f_H (m_l) -f_H(v')) \succ $$
$$\succ (f_H(v), f_H (m_1), ..., f_H(m_l)) =: \mathbf{x}.$$

If one denotes with $\mathbf{z}$ the vector of (unchanged) weights
of groups subordinated to all other non-terminal vertices of $H$,
then, by Lemma \ref{lemma_Zhang_xy},
$\mathbf{f}(T')=(\mathbf{y},\mathbf{z})\succ
(\mathbf{x},\mathbf{z})=\mathbf{f}(H)$.

Assume now that $v \notin g_H(u)$. Then there are disjoint paths $(u, m_1, ..., m_k, m)$ and $(v, m_1', ..., m_l', m)$
(where $k,l\ge0$) in $H$ from vertices $u$ and $v$ to some vertex $m \in M(H)$.

If $f_H(u)> f_H(v)$, then, applying repeatedly formula
(\ref{eq_weights_monotone}) from Lemma
\ref{lemma_Huffman_monotonicity}, we write $f_H(m_i) > f_H(m_i'), i
= 1, ..., \min[k, l]$. It also follows from
(\ref{eq_weights_monotone}) that $k \le l$, since otherwise
$f_H(m_{l+1})> f_H(m)$, which is impossible, since $m_{l+1} \in
g_H(m)$.

Consider a directed tree $T'\in\mathcal{WD}(\mu, d')$ obtained from $H$ by deleting the arc $v'v$ and
adding the arc $v'u$ instead. In the tree $T'$ weights of the groups subordinated to the vertices $u, m_1, ..., m_k$
increase by $f_H(v')$ (i.e., $f_{T'}(u) = f_H(u) + f_H(v')$, $f_{T'}(m_i) = f_H(m_i) + f_H(v'), i = 1, ..., k$),
weights of the groups subordinated to the vertices $v, m_1', ..., m_l'$ decrease by $f_H(v')$ (i.e., $f_{T'}(u) =
f_H(u) - f_H(v')$, $f_{T'}(m_i') = f_H(m_i') - f_H(v'), i = 1, ..., l$), weights of all other vertices (including $m$)
do not change. Therefore, by Lemma \ref{lemma_Zhang_b},
$$\mathbf{y} := (f_{T'}(u), f_{T'}(m_1), ..., f_{T'}(m_k), f_{T'}(v), f_{T'}(m_1'), ..., f_{T'}(m_l'))=$$
$$=(f_H(u)+f_H(v'), f_H(m_1)+f_H(v'), ..., f_H(m_k)+f_H(v'),$$
$$f_H(v)-f_H(v'), f_H(m_1')-f_H(v'), ..., f_H (m_l')-f_H(v'))\succ$$
$$\succ (f_H (u), f_H (m_1), ..., f_H (m_k), f_H (v), f_H(m_1'), ..., f_H (m_l')) =: \mathbf{x}.$$

If $\mathbf{z}$ is a vector of (unchanged) weights of groups
subordinated to all other non-terminal vertices of $H$, then, by
Lemma \ref{lemma_Zhang_xy},
$\mathbf{f}(T')=(\mathbf{y},\mathbf{z})\succ
(\mathbf{x},\mathbf{z})=\mathbf{f}(H)$.

By construction of the Huffman tree, the situation of $f_H(u) = f_H(v)$ is possible only when $d(u)=d(v)$ and
$\mu(u)=\mu(v)$. In this case we cannot use formula (\ref{eq_weights_monotone}) to compare subordinate groups' weights of
elements of both chains, since all possible alternatives of $k = 0$, or $l = 0$, or any sign of the expression
$f_H(m_1) - f_H(m_1')$ in case of $k, l \ge 1$ are possible.

On the other hand, if $f_H(m_1)> f_H(m_1')$, then formula (\ref{eq_weights_monotone}) can be used to show that $k \le
l$, $f_H(m_i) > f_H(m_i')$, $i = 2, ..., k$. In case of the opposite inequality, $f_H(m_1) < f_H(m_1')$, formula
(\ref{eq_weights_monotone}) says that, by contrast, $k \ge l$, $f_H(m_i) < f_H(m_i')$, $i = 2, ..., l$. Repeating this
argument through the chain, we see that only two alternatives are possible:
\begin{itemize}
    \item $0\le p \le k \le l$, $f_H(m_i) = f_H(m_i')$, $i = 1, ..., p$,
    $f_H(m_i) > f_H(m_i')$, $i = p + 1, ..., k$. In this case, as above, we can show that for the directed tree
    $T'\in\mathcal{WD}(\mu, d')$ obtained from $H$ by deleting the arc $v'v$ and adding the arc $v'u$ instead,
    $\mathbf{f}(T')\succ\mathbf{f}(H)$.
    \item $0\le p \le l \le k$, $f_H(m_i) = f_H(m_i')$, $i = 1, ..., p$,
    $f_H(m_i) < f_H(m_i')$, $i = p + 1, ..., l$. In this case the same inequality is true for the directed tree
    $T'\in\mathcal{WD}(\mu, d')$ obtained from $H$ by redirecting to/from vertex $v$ all arcs incident to $u$,
    and by redirecting to/from vertex $u$ all arcs incident to $v$ except the arc $v'v$.
\end{itemize}

Therefore, we proved that a tree $T'\in \mathcal{WD}(\mu, d')$
exists such that $\mathbf{f}(T')\succ \mathbf{f}(H)$. As shown
above, $\mathbf{f}(H)\succeq \mathbf{f}(T)$, so, finally,
$\mathbf{f}(T')\succ \mathbf{f}(T)$.
\end{proof}\end{lemma}

\begin{definition}A directed tree $T \in \mathcal{WD}(\mu, d)$ with a terminal vertex $t$ is called a \emph{proper tree}
if for all $m\in M(T)$, $m\neq t$, $f_T(m) \le \bar{\mu}/2$.\end{definition}

\begin{lemma}\label{lemma_concave} Let a function $\chi(x)$ be concave and increasing for $x\in [0,\bar{\mu}/2]$.
Consider a pair of generating tuples, $\langle\mu,d\rangle$ and $\langle\mu,d'\rangle$, satisfying conditions of Lemma \emph{\ref{lemma_major}}.
If $T\in \mathcal{WD}(\mu,d)$ and $T'\in \mathcal{WD}(\mu,d')$ are directed Huffman trees, then
$$\sum_{v\in V(T')\backslash\{t'\}}\chi(f_{T'}(v))<\sum_{v\in V(T)\backslash\{t\}}\chi(f_T(v)),$$
where $t\in M(T)$ and $t'\in M(T')$ are terminal vertices of $T$ and $T'$ respectively.
\begin{proof}
From Lemma \ref{lemma_major}, such a tree $T''\in
\mathcal{WD}(\mu,d')$ exists that $\mathbf{f}(T'')\succ
\mathbf{f}(T)$. Theorem 2 from \cite{Goubko15} says that
$\mathbf{f}(T')\succeq \mathbf{f}(T'')$, therefore,
$\mathbf{f}(T')\succ \mathbf{f}(T)$. Denote for short $n_1=|V(T)|$,
$\mathbf{f}(T)=\mathbf{f}:=(f_1,...,f_{n_1-1})$,
$\mathbf{f}(T')=\mathbf{f}':=(f'_1,...,f'_{n_1-1})$.

It is known (see Lemma 19 in \cite{Goubko15}) that each directed Huffman tree with degree-monotone weights is a proper tree,
so, $f_T(w)\le \bar{\mu}/2$, $w\in M(\mu,d)\backslash\{t\}$, and $f_{T'}(w')\le \bar{\mu}/2$, $w'\in M(\mu,d')\backslash\{t'\}$.
If a vertex $w\in V(T)$ exists, such that $\mu(w)>\bar{\mu}/2$ (there can be at most one such vertex in $V(T)$),
then $w$ cannot be an internal vertex in $T$ and a pendent vertex in $T'$, since then conditions of Lemma \ref{lemma_major}
imply that $w=v$ and $\mu(u)\ge\mu(w)$, which is impossible. Therefore, $w$ is either a terminal vertex both in $T$
and in $T'$, or a pendent vertex both in $T$ and in $T'$. In the latter case $f_T(v)=f_{T'}(v)=\mu(v)$.

Consequently, $f_i,f_i'\le\bar{\mu}/2$ for $i=1,...,n_1-2$, and if $f_{n_1-1}>\bar{\mu}/2$, then $f_{n_1-1}=f'_{n_1-1}=\mu(w)$.

If $f_{n_1-1}\le \bar{\mu}/2$, the statement of the lemma
follows from Lemma \ref{lemma_Zhang_concave}.

If $f_{n_1-1}> \bar{\mu}/2$, we can write
$$\sum_{v\in V(T)\backslash\{t\}}\chi(f_T(v))-\sum_{v\in V(T')\backslash\{t'\}}\chi(f_{T'}(v))=$$
$$=\sum_{i=1}^{n_1-2}\chi(f_i)+\chi(f_{n_1-1})-\sum_{i=1}^{n_1-2}\chi(f'_i)-\chi(f'_{n_1-1})=\sum_{i=1}^{n_1-2}\chi(f_i)-\sum_{i=1}^{n_1-2}\chi(f'_i).$$

Since $\mathbf{f}'\succ \mathbf{f}$ and $f_{n_1-1}=f'_{n_1-1}$, we
have $(f'_1,...,f'_{n_1-2})\succ (f_1,...,f_{n_1-2})$, and the
statement of the lemma again follows from Lemma
\ref{lemma_Zhang_concave}.
\end{proof}
\end{lemma}

\begin{corollary}\label{corollary_Wiener_1}
If generating tuples $\langle\mu,d\rangle$ and $\langle\mu,d'\rangle$ satisfy conditions of Lemma \emph{\ref{lemma_major}},
then $VWWI^*_{\mathcal{WT}}(\mu,d')<VWWI^*_{\mathcal{WT}}(\mu,d)$.
\begin{proof}
Theorem 3 from \cite{Goubko15} says that $\mathcal{WD}_{VWWI}^*(\mu,d)$ consists of all directed Huffman trees.
Consider the directed Huffman trees $T\in \mathcal{WD}_{VWWI}^*(\mu,d)$, $T'\in \mathcal{WD}_{VWWI}^*(\mu,d')$.
Function $\chi(x)=x(\bar{\mu}-x)$ in (\ref{eq_VWWI_directed})
satisfies the conditions of Lemma \ref{lemma_concave}, so, from (\ref{eq_VWWI_directed}), $VWWI(T')<VWWI(T)$.
Since every proper directed tree from $\mathcal{WD}(\mu,d)$ has a corresponding tree from $\mathcal{WT}(\mu,d)$, and vice versa,
we have $VWWI^*_{\mathcal{WT}}(\mu,d)=VWWI^*_{\mathcal{WD}}(\mu,d)$, and the corollary follows immediately.
\end{proof}
\end{corollary}

In the rest of the section we consider a set $V$ consisting of $n$ vertices with weights $\mu(v)$, $v\in V$. The set $V$ can be thought of as a fixed collection of
(heterogeneous) atoms used as building blocks for molecules. All molecules constructed from these building blocks belong to $\mathcal{WT}(\mu)$.

We want to use Corollary \ref{corollary_Wiener_1}
to show that, similar to Theorem \ref{theorem_C_sufficient}, the vertex-weighted Wiener index is minimized by a tree
having as many vertices of the maximum degree as possible. We cannot apply Corollary \ref{corollary_Wiener_1}
to the whole collection $\mathcal{WT}(\mu)$ of trees with vertices having fixed weights $\mu(\cdot)$
(unless $\mu(v)\equiv const$), as it inevitably contains trees generated by the
 tuples with non-degree-monotone weights, for which Corollary \ref{corollary_Wiener_1} is inapplicable.
 Therefore, we have to carefully limit a set of admissible trees.

\begin{definition}\label{def_extremal_tuple} \sloppy Consider
an admissible set $\mathcal{M}\subseteq \mathcal{WT}(\mu)$ and
denote with $L:=\bigcap_{T\in\mathcal{M}}W(T)$ the set of vertices, which are pendent in all trees from $\mathcal{M}$.
A generating tuple $\langle\mu,\bar{d}\rangle$ is called \emph{extremal} for $\mathcal{M}$, if  $M(\mu, \bar{d})$ consists
of $\lceil\frac{n-2}{3}\rceil$ vertices having the highest weights in $V\backslash L$
(i.e., if $u\in M(\mu, \bar{d})$ and $v\in V\backslash L$, then $\mu(u)\ge\mu(v)$),
at most one vertex $u\in M(\mu, \bar{d})$ has degree $\bar{d}(u)<4$ while others having degree $4$, and, when exists,
$u$ has the minimal weight in $M(\mu, \bar{d})$.
\end{definition}

There can be several extremal tuples for an admissible set, if for some extremal tuple $\langle\mu,\bar{d}\rangle$
such vertices $w\in W(\mu,\bar{d})\backslash L$ and $m\in M(\mu,\bar{d})$
exist that $\mu(w)=\mu(m)$ (swapping $w$ and $m$ then makes a new extremal tuple).
It is clear that weights are degree-monotone in $\langle\mu,\bar{d}\rangle$,
and there are only extremely branched trees in $\mathcal{WT}(\mu,\bar{d})$.

\begin{theorem}\label{theorem_degree_domination}
If $\langle\mu,\bar{d}\rangle$ is an extremal tuple for an
admissible set $\mathcal{M} \subseteq \mathcal{WT}(\mu)$ of
\textbf{chemical} trees with degree-monotone weights, and $H\in
\mathcal{WT}(\mu,\bar{d})$ is a Huffman tree, then $VWWI(H)\le
VWWI_\mathcal{M}^*$, with equality if and only if $\mathcal{M}$
contains $H$ or any other Huffman tree generated by the extremal
tuple.
\begin{proof}
Consider a vertex-weighted tree $T\in \mathcal{M}$. Weights are degree-monotone in $T$,
so, if $d_T(u)>d_T(v)$ and $\bar{d}(u)<\bar{d}(v)$ for some $u,v\in M(T)$, then $\mu(u)=\mu(v)$. Consequently, such
a tree $T'\in \mathcal{WT}(\mu)$ exists that $VWWI(T)=VWWI(T')$, $W(T')=W(T)$, and from $u,v\in M(T)$ and $d_{T'}(u)<d_{T'}(v)$
it follows that $\bar{d}(u)\le \bar{d}(v)$
($T'$ is constructed from $T$ with a permutation of several vertices of equal weight, which does not affect the index value). So, degrees
in the tree $T'$ are ``compatible'' to those in the tuple $\langle\mu,\bar{d}\rangle$.

Assume $T'\in \mathcal{WT}(\mu,d^0)$ for some generating tuple $\langle\mu,d^0\rangle$, and $d^0(\cdot)\not\equiv\bar{d}(\cdot)$.
Construct such a sequence $\langle\mu,d^0\rangle, \langle\mu,d^1\rangle, ..., \langle\mu,d^k\rangle$
of generating tuples with degree-monotone weig\-hts that $k\ge 1$,
$d^k(\cdot)\equiv\bar{d}(\cdot)$, and each pair
$\langle\mu,d^i\rangle,\langle\mu,d^{i+1}\rangle$ of sequential tuples meets the requirements of Lemma~\ref{lemma_major}, $i=0,...,k-1$.

\sloppy We build the elements of this sequence one by one. For a tuple $\langle\mu,d^i\rangle$ define $M^+(\mu,d^i):=\{w\in V: d^i(w)<\bar{d}(w)\}$
and $M^-(\mu,d^i):=\{w\in V: d^i(w)>\bar{d}(w)\}$. Set $d^{i+1}(u^i)=d^i(u^i)+1$ for a vertex
$u^i$ having the highest weight in $M^+(\mu,d^i)$, and $d^{i+1}(v^i)=d^i(v^i)-1$
for a vertex $v^i$ having the least weight in $M^-(\mu,d^i)$, while keeping the degrees of all other vertices.

%
%

%

\begin{figure}[htpb]
\begin{center}
  \includegraphics[width=11cm]{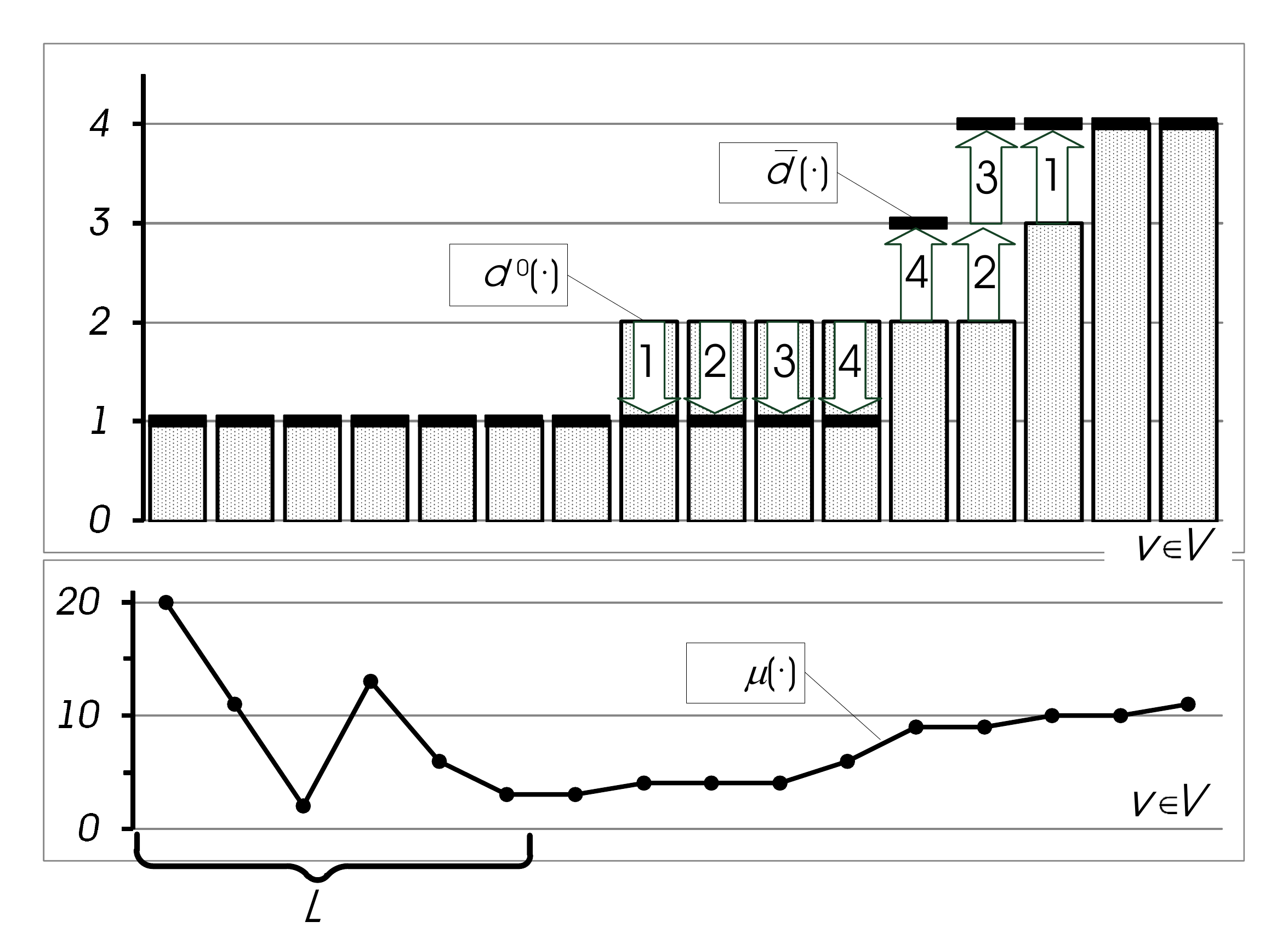}
  \caption{Sequence of degree functions converging to $\bar{d}(\cdot)$. A brace shows $L$.}
  \label{fig_degree_seq}
\end{center}
\end{figure}

The proof is clear from Fig.~\ref{fig_degree_seq}, where numbers on arrows show steps of transformation of the initial degree function.
By Corollary~\ref{corollary_Wiener_1}, $VWWI^*_{\mathcal{WT}}(\mu,\bar{d})<...<VWWI^*_{\mathcal{WT}}(\mu,d^0)\le VWWI(T')=VWWI(T)$.
From Theorem~\ref{theorem_huffman_vwwi} we know that if $T\in \mathcal{WT}(\mu,\bar{d})$, then $VWWI(T)=VWWI^*_{\mathcal{WT}}(\mu,\bar{d})$
if and only if $T$ is a Huffman tree. Therefore, for every $T\in \mathcal{M}$ we have $VWWI^*_{\mathcal{WT}}(\mu,\bar{d})\le VWWI(T)$ with equality
if and only if $T$ is a Huffman tree for $\langle\mu,\bar{d}\rangle$ or another extremal tuple, and the proof is complete.
\end{proof}
\end{theorem}

\begin{figure}[htpb]
\begin{center}
  \includegraphics[width=12cm]{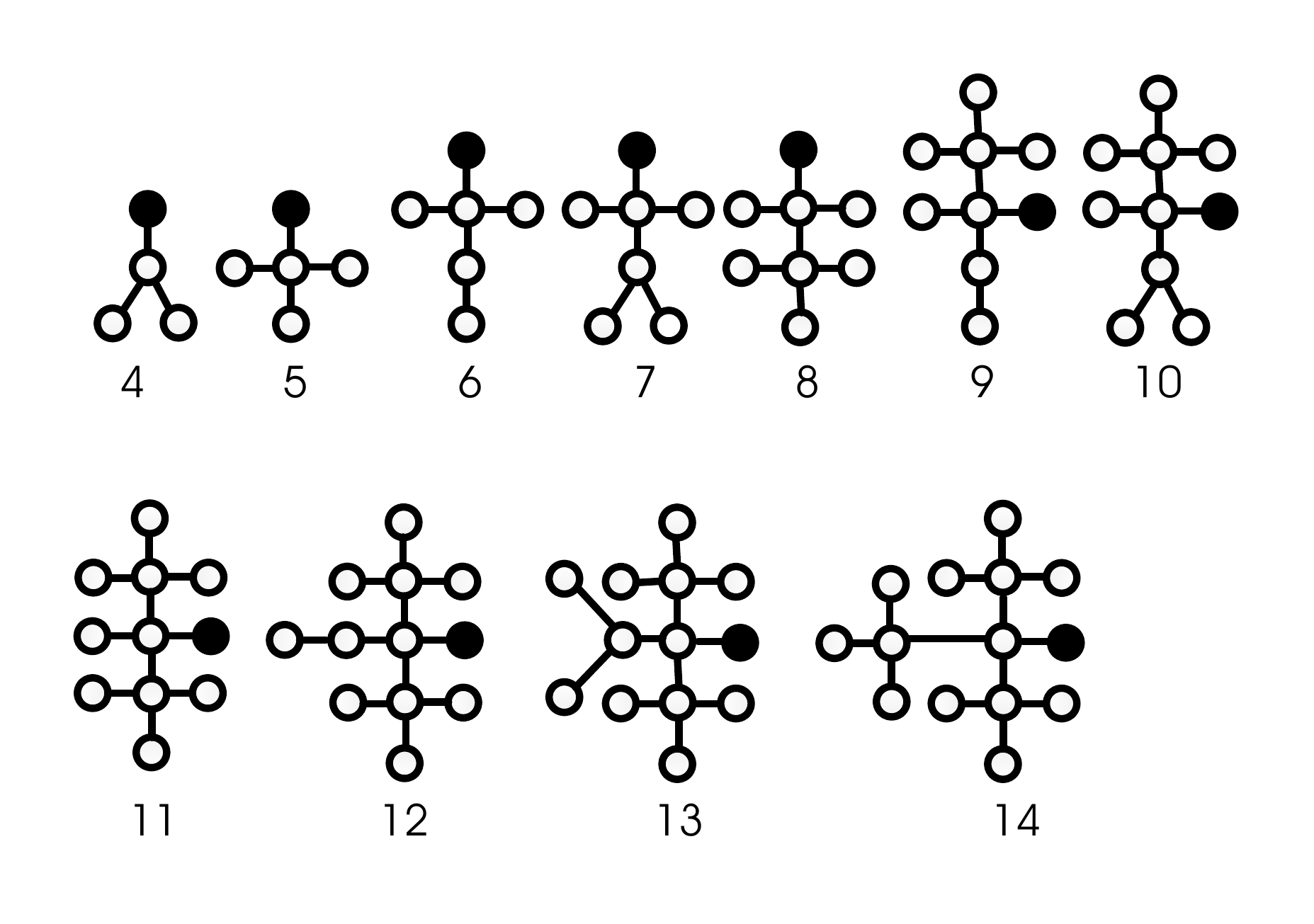}
  \caption{Rooted chemical trees minimizing $WI_\textrm{O}(\cdot)$ for $n=4,...,14$
  (alcohol molecules minimizing $BP^\textrm{II}(\cdot)$ with oxygens filled with black)}
  \label{fig_VWWI_minimizers}
\end{center}
\end{figure}

\begin{corollary}\label{corollary_VWWI_max_degree}
Trees, where $WI_\textrm{\emph{O}}(\cdot)$ achieves its minimum over the set of pendent-rooted trees of order $n=4,...,14$, are depicted in Fig.~\emph{\ref{fig_VWWI_minimizers}}.
\begin{proof}
As we already argued, $WI_\textrm{O}(\cdot)$ can be seen a special case
 of $VWWI(\cdot)$ for the vertex set where all vertices have sufficiently small weight $\varepsilon$ except a
 root having weight $1/\varepsilon$. Let $r$ be the root in all considered pendent-rooted trees.
 Then the set of pendent-rooted trees satisfies conditions of Theorem \ref{corollary_Wiener_1} with $L=\{r\}$, and the
statement of the corollary follows from Theorem \ref{corollary_Wiener_1} and the fact that for a given $n$
the value of the index cannot be improved only for the extremal degree function $\bar{d}(\cdot)$, the one distinct from 1 and 4 at no more than
one vertex. The concrete degree function for each $n$ is justified from (\ref{eq_tree_degrees}), as in Theorem
\ref{theorem_C_sufficient}. From \cite{Goubko15} we know that, for a given degree function $\bar{d}(\cdot)$,
 Huffman trees (shown in Fig.~\ref{fig_VWWI_minimizers}) minimize $VWWI(\cdot)$.
\end{proof}
\end{corollary}


The extension of the above results to trees with vertex degrees limited to any $\Delta\ge 3$ is straightforward.

\section{Minimal Boiling Point}
Below we combine the theorems from the previous section to find an alcohol molecule with the minimal (predicted) boiling point.
The main tool will be the following obvious proposition.

\begin{proposition}\label{prop_subset}
For an arbitrary admissible set $\mathcal{M}$ of graphs and a pair of indices $I(G)$ and $J(G)$ defined for
all graphs $G\in \mathcal{M}$, if $\mathcal{M}_I^*\cap \mathcal{M}_J^*\neq \emptyset$, then $\mathcal{M}_{I+J}^* = \mathcal{M}_I^*\cap \mathcal{M}_J^*$.
\begin{proof} is straightforward, since $\mathcal{M}_I^*\cap \mathcal{M}_J^*$ is a collection of graphs where indices $I(\cdot)$ and $J(\cdot)$
attain their minima simultaneously.
\end{proof}
\end{proposition}

We start with the simplest Regression I (see Equation (\ref{eq_bpI})) being a linear combination of degree-based indices -- the ``ad-hoc index'' $C^\textrm{I}(\cdot)$
 defined with Equation (\ref{eq_ad_hoc_index}) and $S_2(\cdot)$ index,
which penalizes presence of a sub-root of degree 2 in a pendent-rooted tree.

\begin{proposition}\label{prop_BPI}
The collection $\mathcal{R}_{BP^\textrm{\emph{I}}}^*(n)$ of pendent-rooted trees \emph{(}alcohol molecules\emph{)} with the minimal
predicted boiling point $BP^\textrm{\emph{I}}$ for $n = 4,..., 14$ consists
of the trees depicted in Fig.~\emph{\ref{fig_C_minimizers}}
with black circles representing a possible root \emph{(}oxygen atom\emph{)} positions in a tree \emph{(}alcohol molecule\emph{)}.
\begin{proof}
Using Table \ref{tab_param} one easily checks conditions (\ref{eq_degree_index_condition22})-(\ref{eq_degree_index_condition23bis})
of Theorem \ref{theorem_C_sufficient} to hold for the ``ad-hoc index''
$C^\textrm{I}(\cdot)$ and $n\le14$. Therefore, by Corollary \ref{corollary_ad_hoc}, trees from $\mathcal{T}_{C^\textrm{I}}^*(n)$ are
depicted in Fig. \ref{fig_C_minimizers}. We know that the ``ad-hoc'' index does not care for the root position in a tree, and,
therefore, $\mathcal{R}_{C^\textrm{I}}^*(n)$ can be obtained from $\mathcal{T}_{C^\textrm{I}}^*(n)$ by assigning the root to
any pendent vertex of all optimal trees.

From Lemma \ref{lemma_S_i}, $S_2(\cdot)$ is minimized with pendent-rooted trees whose sub-root has degree 3 or 4. Each tree in Fig.~\ref{fig_C_minimizers} has
a pendent vertex being incident to a vertex of degree 3 or 4, so, $\mathcal{R}_{C^\textrm{I}}^*(n)\cap \mathcal{R}_{S_2}^*(n)\neq \emptyset$, and
 conditions of Proposition \ref{prop_subset} hold for all $n=4, ..., 14$. Therefore, $\mathcal{R}_{BP^\textrm{I}}^*(n)$ consists of the trees
 from Fig.~\ref{fig_C_minimizers} with the root assigned to a pendent vertex incident to a vertex of degree 3 or 4. Possible roots are filled with black color in Fig. \ref{fig_C_minimizers}.
\end{proof}
\end{proposition}

The similar reasoning can be carried out for Regression II, which is
defined with Equation (\ref{eq_bpII}) and adds up from the
generalized first Zagreb index $C_1^\textrm{II}(\cdot)$, the cube
root of the distance $WI_\textrm{O}(\cdot)^{\frac{1}{3}}$ of the
oxygen atom in an alcohol molecule, and, again, from $S_2(\cdot)$
index.

\begin{proposition}\label{prop_BPII}
The collection $\mathcal{R}_{BP^\textrm{\emph{II}}}^*(n)$ of rooted trees \emph{(}alcohol molecules\emph{)} with the minimal
predicted boiling point $BP^\textrm{\emph{II}}$ for $n = 4,..., 14$ consists
of the trees depicted in Fig.~\emph{\ref{fig_VWWI_minimizers}}
with a black circle representing a root \emph{(}oxygen atom\emph{)} position in a tree \emph{(}alcohol molecule\emph{)}.
\begin{proof}
Using Table \ref{tab_param} we check that conditions (\ref{eq_degree_index_condition23})-(\ref{eq_degree_index_condition33}) of
Theorem \ref{theorem_C_sufficient} hold, therefore, from Corollary~\ref{corollary_C1}, the
generalized first Zagreb index $C_1^\textrm{II}(\cdot)$ is minimized with \textbf{any} extremely branched tree.

According to Corollary \ref{corollary_VWWI_max_degree},
$\mathcal{R}_{WI_\textrm{O}}^*(n)$ includes only extremely branched
trees depicted in Fig.~\ref{fig_VWWI_minimizers}. Since the cube
root is a monotone function, by Proposition \ref{prop_subset},
$\mathcal{R}_{WI_\textrm{O}^{1/3}+C_1^\textrm{II}}^*(n)=\mathcal{R}_{WI_\textrm{O}}^*(n)$.

The sub-root of all trees depicted in Fig.~\ref{fig_VWWI_minimizers} has degree 3 or 4. Therefore,
$\mathcal{R}_{WI_\textrm{O}}^*(n)\cap \mathcal{R}_{S_2}^*(n)=\mathcal{R}_{WI_\textrm{O}}^*(n)$, and, finally, by Proposition \ref{prop_subset},
$\mathcal{R}_{BP^\textrm{II}}^*(n)$ consists of rooted trees
 from Fig.~\ref{fig_VWWI_minimizers} with roots filled black.
\end{proof}
\end{proposition}

The Basic regression (see Equation (\ref{eq_bp0})) combines the
``ad-hoc'' index $C^\textrm{0}(\cdot)$, the cube root of the
oxygen's distance $WI_\textrm{O}(\cdot)^{\frac{1}{3}}$, and
$S_2(\cdot)$ index. From Corollary \ref{corollary_VWWI_max_degree}
we know that $WI_\textrm{O}(\cdot)$ is minimized with an extremely
branched tree, and so is $WI_\textrm{O}(\cdot)^{\frac{1}{3}}$. At
the same time, neither of the inequalities
(\ref{eq_degree_index_condition23})-(\ref{eq_degree_index_condition23bis})
from Theorem \ref{theorem_C_sufficient} holds for $C^0(\cdot)$, and
we cannot be sure that $C^\textrm{0}(\cdot)$ is minimized with an
extremely branched tree. Therefore, we cannot prove formally that
$BP^0(\cdot)$ achieves its minimum at some extremely branched tree.

Nevertheless, since for both simplified regressions ($BP^\textrm{I}$ and $BP^\textrm{II}$)
an extremely branched tree appears to be optimal, optimality of an extremely branched tree
for the Basic regression $BP^0(\cdot)$ is believed to be a credible hypothesis, which we state below formally.

\begin{conjecture}\label{conjecture_bp0}
If $BP^0(T)=BP^{0*}_\mathcal{R}(n)$ for some pendent-rooted tree $T\in \mathcal{R}(n)$, then $T$ is an extremely branched tree.
\end{conjecture}

\begin{figure}[htpb]
\begin{center}
  \includegraphics[width=11cm]{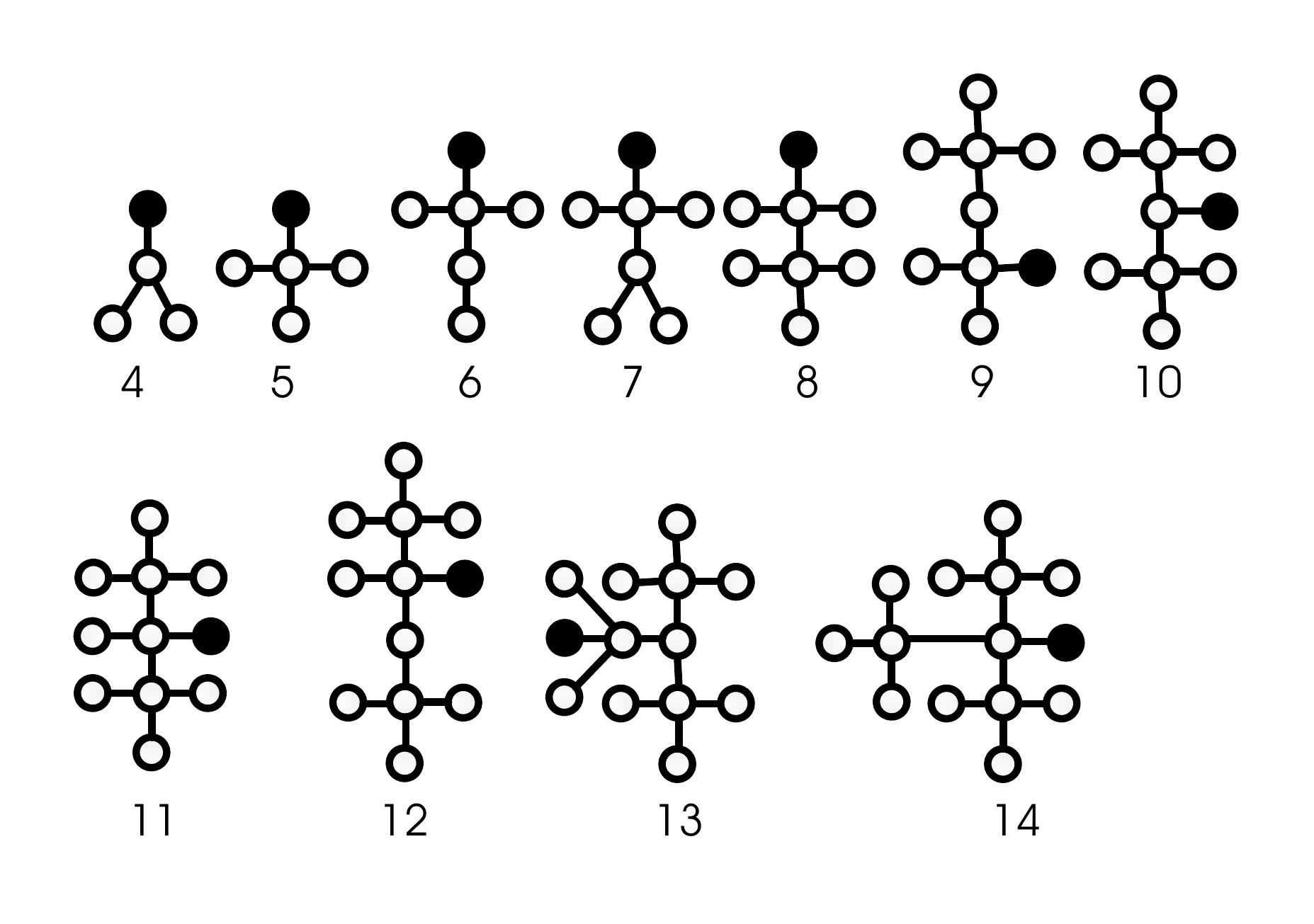}
  \caption{Pendent-rooted chemical trees (alcohol molecules) minimizing $BP^0(\cdot)$ for $n=4,...,14$ with black circles being roots (oxygen atoms)}
  \label{fig_BP0_minimizers}
\end{center}
\end{figure}

\begin{figure}[htpb]
\begin{center}
  \includegraphics[width=14cm]{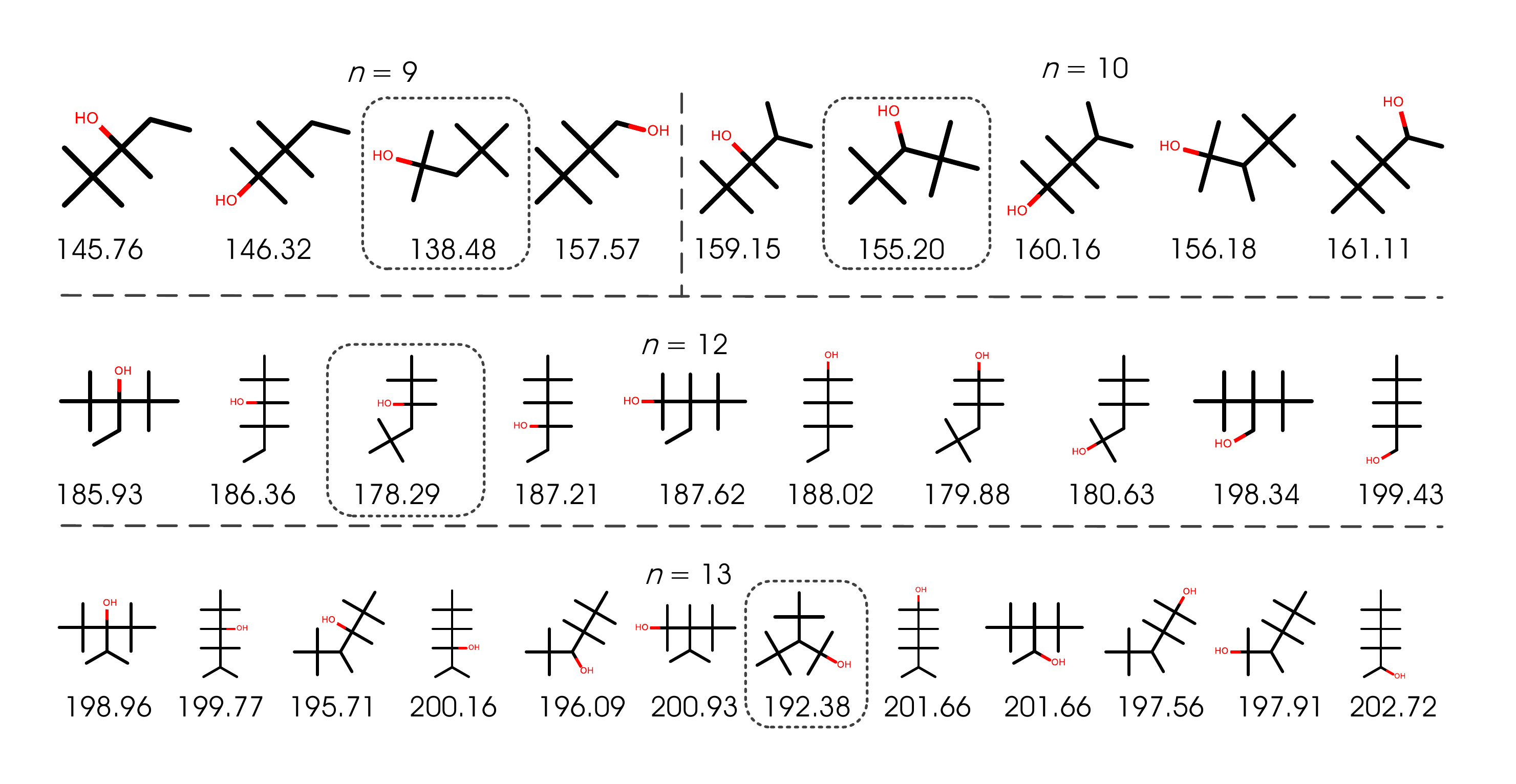}
  \caption{Extremely branched pendent-rooted trees of order $n=9,10,12,13$ and their $BP^0$ values (minima are framed)}
  \label{fig_admissible_sets}
\end{center}
\end{figure}

\begin{proposition}\label{prop_BP0}\sloppy
If Conjecture~\emph{\ref{conjecture_bp0}} is supposed to hold, the collection $\mathcal{R}_{BP^0}^*(n)$ of rooted trees \emph{(}alcohol molecules\emph{)} with the minimal
predicted boiling point $BP^0$ for $n = 4,...,14$ consists
of the trees depicted in Fig.~\emph{\ref{fig_BP0_minimizers}}
with a black circle representing a root \emph{(}oxygen atom\emph{)} position in a tree \emph{(}alcohol molecule\emph{)}.
\begin{proof}
From Conjecture~\ref{conjecture_bp0}, $BP^0$-minimizer is an extremely branched tree.
From Fig.~\ref{fig_C_minimizers} and \ref{fig_VWWI_minimizers} we see that for $n=4,5,6,8$ optimal trees for $WI_\textrm{O}(\cdot)$
and for $C^0(\cdot)$ coincide, while for $n=7,11,14$
$\mathcal{R}_{WI_\textrm{O}}^*(n)\subset\mathcal{R}_{C^0}^*(n)$. Since all sub-roots in $\mathcal{R}_{WI_\textrm{O}}^*(n)$ have
degree 3 or 4, from Proposition~\ref{prop_subset} we conclude that $\mathcal{R}_{BP^0}^*(n)=\mathcal{R}_{WI_\textrm{O}}^*(n)$ for such $n$.

For $n=9,10,12,13$, if Conjecture~\ref{conjecture_bp0} is assumed to hold, both $WI_\textrm{O}$-minimizers and $C^0$-minimizers
are extremely branched trees depicted in Fig.~\ref{fig_C_minimizers} and \ref{fig_VWWI_minimizers}, but $\mathcal{R}_{WI_\textrm{O}}^*(n)\cap
\mathcal{R}_{C^0}^*(n)=\emptyset$, and we cannot use Proposition \ref{prop_subset}. At the same
time, the sets of extremely branched rooted trees of order $n$ are remarkably small for
$n=9,10,12,13$, and are completely enumerated in Fig.~\ref{fig_admissible_sets} along with their predicted normal boiling points
$BP^0$. The tree with the least boiling point for each $n$ is framed in Fig.~\ref{fig_admissible_sets}.
Combining the above findings we obtain Fig.~\ref{fig_BP0_minimizers}.
\end{proof}
\end{proposition}

Predictions of the minimal boiling point based on the Wiener index slightly differ from those
based on the ``ad-hoc'' index $C^\textrm{I}(\cdot)$.
Minimizers of both indices are extremely branched trees. For $n=4,5,6,8$ optimal trees for $WI_\textrm{O}(\cdot)$ and for $C^\textrm{I}(\cdot)$ coincide. For $n=7,11,14$
$\mathcal{R}_{WI_\textrm{O}}^*(n)\subset \mathcal{R}_{C^\textrm{I}}^*(n)$, so, $WI_\textrm{O}(\cdot)$ appears to be a refinement of $C^\textrm{I}(\cdot)$,
and Regression II refines Regression I.

For $n=9,10,12,13$ minimizers of $WI_\textrm{O}(\cdot)$ and of $C^\textrm{I}(\cdot)$ differ. By Corollaries \ref{corollary_ad_hoc} and \ref{corollary_VWWI_max_degree},
extremely branched trees minimize both indices. They have the only internal vertex $v$ of non-maximal degree, but $C^\textrm{I}(\cdot)$ suggests
settling this vertex in the very center of a molecule surrounding it with the other internal vertices, while $WI_\textrm{O}(\cdot)$
says vertex $v$ must be a \emph{stem vertex} incident to only one internal vertex. Therefore, predictions of Regression II differ from those of Regression I.

Basic regression, which contains a weighted linear combination of the Wiener index and the ``ad-hoc'' index, represents a sort of
intermediate behavior between these extremal trends (at least for extremely branched trees under Conjecture~\ref{conjecture_bp0}).
The weight of the Wiener index in the regression is not enough to move the
minimal tree sufficiently from the trees depicted in Fig.~\ref{fig_C_minimizers}, and Basic regression becomes a yet
another refinement of Regression~I.


\section{Conclusion}

The focus of this paper is development of optimization techniques
for combinations of some well-known and novel topological indices
over chemically interesting sets of graphs. We derived conditions
under which an extremely branched tree minimizes the sum of the
second Zagreb index and of the generalized first Zagreb index. We
also found minimizers of the vertex-weighted Wiener index over the
set of chemical trees with given vertex weights.

We enumerated index minimizers for moderate (up to 14) non-hydrogen atom count in a molecule,
and combined them in several regressions of different complexity to forecast a simple
alcohol molecule with the lowest boiling point.

For simpler regressions (Regressions I and II) we managed to obtain a complete analytical characterization of extremal alcohol molecules,
while for the most complex (yet the most precise)
``basic'' regression we had to limit our attention to extremely branched trees (see Conjecture~\ref{conjecture_bp0}) and employed
the brute-force enumeration to find molecules of low-boiling alcohols.

Forecasts based on different regressions slightly differ, but they all comply with
the collected experimental data on normal boiling points of simple alcohols.

Finally, let us sketch several promising directions of future
research. An obvious shortcoming of this paper is
Conjecture~\ref{conjecture_bp0}, which is explained but is not
proven formally. To justify it, we need to refine sufficiently our
optimization techniques, \emph{viz}, to optimize jointly the Wiener
index and the Zagreb indices.

On the other hand, there is wide space for investigation of popular
combinations of topological indices forecasting important physical
and chemical properties of compounds.

\end{document}